\documentclass{article}
\usepackage{amssymb}
\newcommand{\Rn}{\mathbf{R}}
\newcommand{\Cn}{\mathbf{C}}
\newcommand{\ord}{\mathrm{ord}\,}

\newcommand{\ini}{\mathrm{in}\,}
\newcommand{\Zer}{\mathrm{Zer}}

\newcommand{\Teis}[2]{
   \setlength{\unitlength}{1ex}
   \begin{picture}(2,0)(0,0.4)
      \put(0,1.1){\line(1,0){2}}
      \put(0,0.9){\line(1,0){2}}
      \put(1,1.2){\makebox(0,0)[b]{$\scriptstyle #1$}}
      \put(1,0.8){\makebox(0,0)[t]{$\scriptstyle #2$}}
   \end{picture}}

\newcommand{\Teisssr}[4]{
   \setlength{\unitlength}{1ex}
   \begin{picture}(#3,3)(0,0.4)
      \put(0,1.15){\line(1,0){#3}}
      \put(0,0.85){\line(1,0){#3}}
      \put(#4,1.3){\makebox(0,0)[b]{$#1$}}
      \put(#4,0.7){\makebox(0,0)[t]{$#2$}}
   \end{picture}}

   \newtheorem{Theorem}{Theorem}
\newtheorem{Lemma}{Lemma}
\newtheorem{Corollary}{Corollary}
\newtheorem{Remark}{Remark}

\newtheorem{Properties}{Properties}
\newtheorem{Property}{Property}

\newtheorem{Proposition}{Proposition}
\newtheorem{Definition}{Definition}
\newtheorem{Example}{Example}

\title{Characterization of jacobian Newton polygons of plane
branches and new criteria of irreducibility \footnotetext{
     \noindent   \begin{minipage}[t]{5in}
       {\small
       2000 {\it Mathematics Subject Classification:\/} Primary 32S55;
       Secondary 14H20.\\
       Key words and phrases: irreducible plane curve,
       jacobian Newton polygon, polar invariant, approximate root.\\
       The first-named author was partially supported by the Spanish Project
       PNMTM 2007-64007.}
       \end{minipage}}}

\author{Evelia R.\ Garc\'{\i}a Barroso and Janusz Gwo\'zdziewicz}

\begin{document}
\maketitle
\begin{center}
{\it Dedicated to Professor Arkadiusz P\l oski on his
60$^{\mbox{\scriptsize th}}$ birthday}
\end{center}

\begin{abstract}
\noindent In this paper we characterize, in two different ways, the
Newton polygons which are jacobian Newton polygons of a branch.
These characterizations give in particular combinatorial criteria of
irreducibility for complex series in two variables and necessary
conditions which a complex curve has to satisfy in order to be the
discriminant of a complex plane branch.
\end{abstract}

\section{Introduction}

\medskip

\noindent Teissier in \cite{Teissier1} introduced the notion of {\em
jacobian Newton polygon} which is the Newton polygon {\em in the
coordinates} $(u,v)$ of the discriminant, i.e. the image of the
critical locus, of a map $$(\ell ,f)\colon (\mathbf
C^{n+1},0)\longrightarrow (\mathbf C^2,0)$$ given by
$(u,v)=(l(x_0,\ldots,x_n),f(x_0,\ldots,x_n))$, where $l$ is a
sufficiently gene\-ral linear form and $f(x_0,\ldots,x_n)$ is a
convergent power series and he proved that these jacobian Newton
polygons are constant for the members of a family of equisingular
germs of complex analytic isolated hypersurface singularities.

\medskip

\noindent The inclinations of the compact edges of these jacobian
Newton polygons are rational numbers called the {\em polar
invariants} of the germ.

\medskip

\noindent For germs of plane complex analytic curves, with the usual
definitions of equisingularity, and in the case of a germ of plane
irreducible curve (i.e. a branch), Merle shows in \cite{Merle} that
the datum of the jacobian Newton polygon determines and is
determined by the equisingularity class of the curve (or
equivalently its embedded topological type). The formulae of Merle
have been generalized to the case of {\em reduced} plane curve germs
by Casas, Eggers, Garc\'{\i}a Barroso, Gwo\'zdziewicz-P\l oski,
Maugendre and Wall among others (see \cite{Casas}, \cite{Eggers},
\cite{GB}, \cite{GP3}, \cite{Maugendre1}, \cite{Maugendre2} and
\cite{Wall}) and they depend only on the equisingularity class of
the curve. In contrast to the usual Newton polygon, the jacobian
Newton polygon is independent of the choice of coordinates, and
encodes a lot of information about the local geometry of a plane
curve (see 4.3 of \cite{Teissier1} for the irreducible case), for
example the \L ojasiewicz exponents for the inequalities $\vert
\hbox{\rm grad} f(z) \vert \geq C_1 \vert z \vert^{\theta}$ and
$\vert \hbox{\rm grad} f(z) \vert \geq C_2 \vert f(z)
\vert^{\theta}$, with $z$ near $0\in \mathbb C^2$ and $C_1$, $C_2$
constants (see \cite{Teissier1}, Corollaire 2, page 270). The
behavior of the curvature of the Milnor fibers is also determined by
the jacobian Newton polygon (see \cite{Garcia-Teissier}).

\medskip
\noindent In 1982, P. Maisonobe in \cite{Maisonobe}, gave necessary
and sufficient conditions on Puiseux exponents of an irreducible
plane curve germ of $\mathbf C^2$ so that it is the discriminant
curve of an analytic map germ $(l,f):\mathbf C^2 \longrightarrow
\mathbf C^2$ where $f$ has an irreducible critical locus. In this
paper we characterize, in two different ways, the special convenient
Newton polygons which are jacobian Newton polygons of a branch,
answering a question of A. Lenarcik and A. P\l oski. These
characterizations give in particular combinatorial criteria of
irreducibility for complex series in two variables and a necessary
condition which a complex curve has to satisfy in order to be the
discriminant of a complex plane branch.

\medskip

\noindent During the preparation of this work the first named author
enjoyed the hospitality of the Technical University of Kielce and
the Institute of Mathematics of Paris.

\medskip

\section{Plane branches, Puiseux expansions and semigroup}
\label{section-Puiseux}
\medskip

\noindent For us, a branch is an irreducible germ of a complex
analytic curve. A plane branch is given by a convergent power series
$f(x,y)\in \mathbf C\{x,y\}$ which is not a unit and is irreducible
in that ring. The branch is the germ at $0$ of the set of solutions
of $f(x,y)=0$.

\medskip

\noindent By the theorem of Newton, after possibly a change of
coordinates to achieve that $x=0$ is transversal to it at $0$, the
branch $C$ can be parametrized near $0$ as follows
$$\begin{array}{llr}x(t)&=t^d\cr
y(t)&=a_et^e+a_{e+1}t^{e+1}+\cdots
+a_{e+j}t^{e+j}+\cdots\hskip1truecm \hbox{with}\ e\geq
d\end{array}$$

\noindent or equivalently by
$$
y(x^{1/d})= a_ex^{e/d}+a_{e+1}x^{(e+1)/d}+\cdots
+a_{e+j}x^{(e+j)/d}+\cdots\hskip1truecm \hbox{with}\ e\geq d$$

\noindent where $d$ is the order of the series $f(x,y)$.

\medskip

\noindent According to Puiseux the branch $C$ admits $d$ different
Newton expansions $\{y_i(x^{1/d})\}_{i=1}^d$, with $y_i(x^{1/d})=
y(\omega_i x^{1/d})$ where $\omega_i$ are the $d$-roots of the unity
in $\mathbf C$. Moreover we may write $f$ as the product
$$ f(x,y)=\mathit{unit}\prod_{i=1}^d(y-y_i(x^{1/d})).$$

\medskip

\noindent The above expansions are called {\em Puiseux roots of the
branch $C$.}

\medskip

\noindent It is well-known that if $d>1$ then there exists $g\in
\mathbf N\setminus\{0\}$ and integer numbers
$\beta_1<\beta_2<\ldots<\beta_g$ such that $\hbox{\rm
g.c.d.}(d,\beta_1,\ldots,\beta_g)=1$ and the orders of the series
$y_i(x^{1/d})-y_j(x^{1/d})$ with $i\neq j$ are exactly the rational
numbers $\left\{\frac{\beta_i}{d}\right\}_{i=1}^g$. The sequence
$(\beta_0=d,\beta_1,\ldots,\beta_g)$ is called the {\em Puiseux
characteristic} or {\em Puiseux exponents of the branch} $C$ defined
by the equation $f(x,y)=0$.

\medskip

\noindent Recall that an increasing sequence $\delta_0<\ldots
<\delta_h$ of positive integer numbers is the Puiseux characteristic
of a branch if and only if $\hbox{\rm
g.c.d.}(\delta_0,\ldots,\delta_h)=1$ and $\hbox{\rm
g.c.d.}(\delta_0,\ldots,\delta_k)< \hbox{\rm
g.c.d.}(\delta_0,\ldots,\delta_{k-1})$ for all $k\in
\{1,\ldots,h\}$.

\medskip

\noindent The {\em semigroup} associated to  the branch $f(x,y)=0$
is by definition $\Gamma:=\{\hbox{\rm
ord}g(t^d,y(t))\;:\;g\not\equiv 0 \;(\hbox{\rm mod} f)\}\subseteq
\mathbf N$.

\medskip

\noindent Zariski (see \cite{Zariski}) proved that the semigroup
$\Gamma$ admits a minimal set of generators $
\overline{\beta}_0<\overline{\beta}_1<\cdots < \overline{\beta}_g$,
that is

$$\Gamma=\bigl\langle
\overline{\beta}_0,\overline{\beta}_1,\ldots
,\overline{\beta}_g\bigr\rangle:=\overline{\beta}_0\mathbf
N+\overline{\beta}_1\mathbf N+\cdots +\overline{\beta}_g\mathbf N.$$

\noindent This set of generators is uniquely determined by the
semigroup $\Gamma$, and determines it.

\medskip

\noindent On the other hand, according to Bresinsky (see Theorem 2,
p. 383 of \cite{Bresinsky}), a semigroup $\Gamma$ of positive
integer numbers generated by $\gamma_0, \ldots,\gamma_r$, is the
semigroup of a plane branch iff it verifies the following three
properties:

\begin{enumerate}
\item $\hbox{\rm g.c.d.}(\gamma_0,\ldots,\gamma_r)=1$,
\item $\hbox{\rm
g.c.d.}(\gamma_0,\ldots,\gamma_i)< \hbox{\rm
g.c.d.}(\gamma_0,\ldots,\gamma_{i-1})$ for all $i\in
\{1,\ldots,r\}$, and
\item $\frac{\hbox{\rm g.c.d.}(\gamma_0,\ldots,\gamma_{i-1})}
{\hbox{\rm g.c.d.}(\gamma_0,\ldots,\gamma_i)}\gamma_i<\gamma_{i+1}$
for $1\leq i \leq r-1$.
\end{enumerate}

\medskip
\noindent Zariski (\cite{Zariski}) proves that the Puiseux
characteristic and the minimal set of the generators of the
semigroup of a branch are equivalent. More precisely if the Puiseux
characteristic of the branch $C$ is $(\beta_0,\ldots,\beta_g)$ and
the generators of its semigroup are
$\overline{\beta}_0,\ldots,\overline{\beta}_g$ then
$\overline{\beta}_i=\beta_i$ for $i\in \{0,1\}$ and
$\overline{\beta}_{i+1}=n_i\overline{\beta}_i+\beta_{i+1}-\beta_i$
for $i\in\{1,\ldots,g-1\}$ where $n_i=\frac{\hbox{\rm
g.c.d}(\beta_0,\ldots,\beta_{i-1})}{\hbox{\rm
g.c.d}(\beta_0,\ldots,\beta_{i})}$. Further, in the proof of
Corollary~\ref{Cor3} we will need an arithmetic property of the
semigroup $\Gamma$:

\begin{Property}
\begin{equation}\label{arith_prop}
\hbox{\rm g.c.d.}\left(\overline{\beta}_0,
(n_1-1)\overline{\beta}_1, \ldots,(n_g-1)\overline{\beta}_g\right)=1
\end{equation}
\end{Property}

\noindent \textbf{Proof.} We prove by induction on~$k$ that
$$\hbox{\rm g.c.d.}\left(\overline{\beta}_0, (n_1-1)\overline{\beta}_1,
\ldots,(n_k-1)\overline{\beta}_k\right) = \hbox{\rm
g.c.d.}\left(\overline{\beta}_0, \overline{\beta}_1, \ldots,
\overline{\beta}_k\right) .$$

\noindent By inductive hypothesis
\begin{eqnarray*}
&& \hbox{\rm g.c.d.}\left(\overline{\beta}_0,
(n_1-1)\overline{\beta}_1,
\ldots, (n_k-1)\overline{\beta}_k, (n_{k+1}-1)\overline{\beta}_{k+1}\right) \\
\quad &=& \hbox{\rm g.c.d.}\left(\hbox{\rm
g.c.d.}(\overline{\beta}_0,(n_1-1)\overline{\beta}_1,
\ldots,(n_k-1)\overline{\beta}_k) , (n_{k+1}-1)\overline{\beta}_{k+1}\right) \\
\quad &=& \hbox{\rm g.c.d.}\left(\hbox{\rm
g.c.d.}(\overline{\beta}_0,\overline{\beta}_1,
\ldots,\overline{\beta}_k) , (n_{k+1}-1)\overline{\beta}_{k+1}\right) \\
\quad &=& \hbox{\rm g.c.d.}\left(\hbox{\rm
g.c.d.}(\overline{\beta}_0,\overline{\beta}_1,
\ldots,\overline{\beta}_k) , -\overline{\beta}_{k+1}\right) =
\hbox{\rm g.c.d.}\left(\overline{\beta}_0,\overline{\beta}_1,
\ldots,\overline{\beta}_{k+1}\right) .
\end{eqnarray*}
In above computation we used elementary properties of the greatest
common divisor and the relation $n_{k+1}\overline{\beta}_{k+1} =
 \frac{\overline{\beta}_{k+1}}{\hbox{\rm g.c.d.}(\overline{\beta}_0,\overline{\beta}_1,
   \ldots,\overline{\beta}_{k+1})}\hbox{\rm g.c.d.}(\overline{\beta}_0,\overline{\beta}_1,
   \ldots,\overline{\beta}_k)$.

\bigskip

\noindent Both Puiseux characteristic and  semigroup form  complete
sets of topological invariants of a branch. It means that they
characterize the equisingularity class of the branch (see 4.3 of
\cite{Teissier1}).

\section{Newton polygons}
\label{Newton polygons}
\medskip

\noindent Let $\mathbf R_+=\{x\in\mathbf R:\,x\geq 0\}$. Any two
subsets $A,B$ of the first quadrant $\mathbf R_+^2$ can be added
coordinate-wise, to give the Minkowski sum  $A+B=\{a+b:\,a\in
A\mbox{ and }b\in B\}$ of $A$ and $B$.
The subset ${\cal N}$ of $\mathbf R_+^2$ is a {\em Newton polygon}
if ${\cal N}$ is the convex hull of  $S+\mathbf R_+^2$ for some
$S\subset\mathbf N^2$, and in this case we denote it by ${\cal
N}(S)$. The boundary of a Newton polygon is a broken line with
infinite horizontal and vertical sides, possibly different from the
coordinate axis and a finite number of compact edges.

\medskip

\noindent According to Teissier (see \cite{Teissier2} and
\cite{Teissier3}) we put $\{\Teis{a}{b}\}={\cal
N}(\{(a,0),(0,b)\})$, $\{\Teis{a}{\infty}\}={\cal N}(\{(a,0)\})$ and
$\{\Teis{\infty}{b}\}={\cal N}(\{(0,b)\})$ for any $a,b>0$ and call
such polygons {\em elementary Newton polygons}. If $\{\Teis{a}{b}\}$
is an elementary Newton polygon, with $a\neq+\infty$, $b\neq
+\infty$, then  its {\em inclination} is by definition the rational
number $\frac{a}{b}$.


\medskip

\noindent The Newton polygons form a semigroup (see
\cite{Teissier2}, Section 3.6, page 616) with the Minkowski sum, and
the elementary Newton polygons generate it. Every Newton polygon
${\cal N}$ can be written, in a unique way, called {\em canonical
form}, as a finite sum

$${\cal N}=\sum_{i=1}^r\left\{\Teis{a_i}{b_i}\right\}$$

\noindent where the inclinations of the terms form an strictly
increasing sequence.

\medskip

\noindent The {\em height} of the Newton polygon ${\cal N}$ is by
definition the length of the projection of compact edges of ${\cal
N}$ on the vertical coordinate axis which we will denote $\hbox{\rm
ht}({\cal N})$. We have $\hbox{\rm ht}({\cal
N})=\sum_{i=1}^r\hbox{\rm ht}(\{\Teis{a_i}{b_i}\})=\sum_{i=1}^r
b_i$.

\medskip

\noindent A Newton polygon is {\em convenient} if intersects both
coordinate axis (see \cite{Kouchnirenko}) and it is {\em special} if
the inclinations of its compact faces are greater than $1$ and
intersects the vertical coordinate axis.

\medskip

\noindent Note that any convenient Newton polygon ${\cal
N}=\sum_{i=1}^r\{\Teis{a_i}{b_i}\}$ is determined by the
inclinations $\left\{\frac{a_i}{b_i}\right\}_{i=1}^r$ of its compact
faces and their respective heights $\{b_i\}_{i=1}^r$.

\medskip

\noindent Fix a complex nonsingular surface i.e. a complex
holomorphic variety of dimension 2. Let $(x,y)$ be a chart centered
at $O$. In all this paper we consider {\em reduced} plane curve
germs, that is curves with local equations  $f(x,y)=\sum a_{ij}x^i
y^j \in\mathbf C\{x,y\}$ without multiple factors. We put ${\cal
N}_{x,y}(C)={\cal N}(S)$ where $S=\{(i,j)\in\mathbf N^2:\,a_{ij}\neq
0\}$ which is called {\em Newton polygon of the curve $C\equiv
f(x,y)=0$ in the coordinates $(x,y)$.} Clearly ${\cal N}_{x,y}(C)$
depends on $(x,y)$.

\medskip

\noindent We can extend the definitions of this section to any
subset $S\subseteq \mathbf Q^2$ such that there exists a positive
integer $m$  with $m.S\subseteq \mathbf N^2$. The Newton polygons
obtained by this gene\-ralization will be called {\em rational
Newton polygons}. The usual Newton polygon with integral vertices
will be called {\em integral Newton polygons.}
\section{Main result}

\medskip

\noindent Let $f=f(x,y)\in\Cn\{x,y\}$ be a power series of order~$d$
without multiple factors such that the vertical axis is transverse
to the curve $f(x,y)=0$. Using Puiseux factorisation we may write
$f$ and its partial derivative $f_y'$ as the products
\begin{equation}\label{Eq1}
   f(x,y)=\mathit{unit}\prod_{i=1}^d(y-\alpha_i(x)),
\end{equation}
$$ f_y'(x,y)=\mathit{unit}\prod_{j=1}^{d-1}(y-\gamma_j(x)), $$

\noindent where $\alpha_i(x)$, $\gamma_j(x)$ are fractional power
series, that is, elements of the ring $\mathbf
C\{x\}^{*}=\bigcup_{n\in \mathbf N}\mathbf C\{x^{1/n}\}.$

\medskip

\noindent If $f'_y=g_1\cdots g_s$ is the factorization of $f'_y$ in
irreducible factors, Teissier proves (\cite{Teissier1}) that the
jacobian Newton polygon ${\cal N}_J(f)$ of $f(x,y)=0$ equals
$$\sum_{j=1}^{s}
\left\{\Teisssr{(f,g_j)_0} {\ord g_j}{8}{4}\right\} ,
$$

\noindent where $(f,g_j)_0$ denotes the intersection number of $f$
and $g_j$, where different elementary Newton polygons may have the
same inclinations.

\medskip

\noindent We can rewrite the above sum in terms of Puiseux roots of
$f'_y$. By Zeuthen's rule (c.f. e.g. \cite{Ploski} Proposition 2.1)
for every $j\in \{1\dots s\}$ we have
$\left\{\Teisssr{(f,g_j)_0}{\ord g_j}{8}{4}\right\}=\sum_i
\left\{\Teisssr{\ord f(x,\gamma_i(x))}{1}{16}{8}\right\}$ where the
sum runs over $\gamma_i(x)$ which are Puiseux roots of $g_j$.

\medskip

\noindent So we can write
$${\cal N}_J(f)=\sum_{j=1}^{d-1}
\left\{\Teisssr{\hbox{\rm ord} f(x,\gamma_j(x))}{1}{14}{7}\right
\}$$

\noindent as a rational Newton polygon.

\medskip

\noindent The polar invariants of $f$, that is, the inclinations of
${\cal N}_J(f)$ are ${\hbox{\rm ord} f(x,\gamma_j(x))}$. If $q_j$ is
the polar invariant of $f$ associated to the edge $\Gamma$ of ${\cal
N}_J(f)$, we call the height of $\Gamma$ the {\em multiplicity} of
$q_j$ and we denote it by $m_j$, consequently ${\cal N}_J(f)$ is
determined by the unordered sequence $Q(f)=\langle\,
q_0:m_0,q_1:m_1,\dots,q_r:m_r\, \rangle$. So the canonical form of
${\cal N}_J(f)$ is
$${\cal N}_J(f)=\sum_{j=0}^{r} \left\{\Teisssr{m_jq_j}{m_j}{8}{4}\right
\}.$$

\medskip

\noindent It is well-known (see \cite{Teissier1}) that ${\cal
N}_J(f)$ is determined by the equisingularity class of $\{f=0\}$. If
$\{f=0\}$ is a branch having semigroup
$\bigl\langle\overline{\beta}_0,\overline{\beta}_1,\ldots,\overline{\beta}_g\bigr\rangle$
then by Merle (see \cite{Merle}) the canonical form of ${\cal
N}_J(f)$ is
\begin{equation}\label{Me}
  {\cal N}_J(f)=\sum_{k=1}^g\left\{\Teisssr{(n_k-1)\overline{\beta}_k}
                                       {(n_k-1)n_1\dots n_{k-1}}{18}{9}
                                       \right\}\;.
\end{equation}
Using~(\ref{Me}) we can compute
$\bigl\langle\overline{\beta}_0,\overline{\beta}_1,\ldots,
\overline{\beta}_g\bigr\rangle$ from ${\cal N}_J(f)$. Hence for $f$
irreducible ${\cal N}_J(f)$ determines the equisingulatity class of
$\{f=0\}$.

\medskip

\noindent Our aim is to show that jacobian Newton polygons of
irreducible series are distinguished among all jacobian Newton
polygons.

\begin{Theorem}\label{Th1}
Let $f,g\in\Cn\{x,y\}$ be such that ${\cal N}_J(f)={\cal N}_J(g)$
and assume that $f$ is irreducible. Then $g$ is also irreducible.
\end{Theorem}

\medskip

\noindent The aim of the next three sections  is to prove Theorem
\ref{Th1}.


\section{The Kuo-Lu Lemma}
\label{Kuo-Lu lemma} \noindent Let $\phi$, $\psi$ be fractional
power series of variable~$x$. We will denote the order of difference
$\phi(x)-\psi(x)$ by $O(\phi, \psi)$ and call it the {\em contact
order.}

\medskip

\noindent For every $\phi_1$, $\phi_2$, $\phi_3\in\Cn\{x\}^{*}$ we
have (see for example page 69 of \cite{Wall2}):

\begin{equation}\label{Eq2}
   \mbox{if $O(\phi_1,\phi_2)\leq O(\phi_2,\phi_3)\leq O(\phi_1,\phi_3)$
         then $O(\phi_1,\phi_2)=O(\phi_2,\phi_3).$}
\end{equation}

\noindent Let $\phi\in\Cn\{x\}^{*}$ and
$r\in\Rn^{+}\cup\{+\infty\}$. We will call the set
$$ B(\phi,r) = \{\,\psi\in\Cn\{x\}^{*}: O(\phi,\psi) \geq r \,\}
$$
a \textit{pseudo-ball with center at $\phi$ and radius $r$.}

\noindent It follows from (\ref{Eq2}) that pseudo-balls have the
following
metric properties:\\[1ex]
(i) every element of a pseudo-ball is its center,\\
(ii) if $B_1=B(\phi_1,r_1)$, $B_2=B(\phi_2,r_2)$ are two disjoint
pseudo-balls
     and $\psi_1\in B_1$, $\psi_2\in B_2$ then $O(\phi_1,\phi_2)=O(\psi_1,
\psi_2)$ i.e.
     contact order between elements of two disjoint pseudo-balls
     does not depend on the choice of these elements,\\
(iii) for all pseudo-balls  $B_1=B(\phi_1,r_1)$ and
$B_2=B(\phi_2,r_2)$
     one of three possibilities holds: \\
$B_1\cap B_2=\emptyset$ \quad $\,$ if $O(\phi_1,\phi_2)<\min(r_1,r_2)$, \\
$B_1\subset B_2$ \quad  \quad  \quad if
$O(\phi_1,\phi_2)\geq\min(r_1,r_2)$ and
$r_1\geq r_2$, \\
$B_2\subset B_1$ \quad \quad  \quad if
$O(\phi_1,\phi_2)\geq\min(r_1,r_2)$ and
$r_2\geq r_1$. \\

\noindent Take a power series $f(x,y)$ without multiple factors of
the form~(\ref{Eq1}), let $\Zer f=\{\alpha_1,\dots,\alpha_d\}$ be
the set of its Puiseux roots and put
$$ {\cal B} := \{\,B(\alpha,O(\alpha,\alpha')): \alpha, \alpha'\in \Zer f\,\},
$$

\noindent which is a partially ordered set with the inclusion
operation.

\noindent We will denote by $h(B)$ the radius of $B\in \cal B$ and
call this number the \textit{height} of $B$.

\medskip

\noindent Inclusion relation gives $\cal B$ a structure of a tree
called \textit{the Kuo-Lu tree-model}~$T(f)$ (see \cite{Kuo-Lu}).
The {\em root} of $T(f)$ is the pseudo-ball of the minimal height
$\min O(\alpha_i,\alpha_j)$ which contains all Puiseux roots of $f$.
{\em Leaves} of $T(f)$ are pseudo-balls
$B(\alpha_i,+\infty)=\{\alpha_i\}$ of infinite heights which can be
identified with Puiseux roots of~$f$. A path from the root to the
leave $\{\alpha_i\}$ connects succesive $B\in \cal B$ of increasing
heights for which $\alpha_i\in B$. Finally a pseudo-ball $B_1$ is a
{\em child} of the pseudo-ball $B_2$ if $B_1\subsetneq B_2$ and
there is not a pseudo-ball $B\in\cal B$ such that $B_1\subsetneq B
\subsetneq B_2$.

\medskip

\noindent Let $B$ be an element of $\cal B$ and let $\gamma$ be any
fractional power series. We will say that $\gamma$ \textit{grows
from} $B$ if and only if $\gamma\in B$ and $O(\gamma,\alpha)=h(B)$
for all $\alpha\in\Zer f\cap B$. In this case contact orders between
$\gamma$ and Puiseux roots of $f$ are given by
$O(\gamma,\alpha)=O(B,\alpha)$ where
$$  O(B,\alpha) :=
\left\{\begin{array}{ll}
    h(B)              & \mbox{if $\alpha\in B$} \\
    O(\alpha',\alpha) & \mbox{otherwise}
\end{array}\right.
$$
and $\alpha'$ is any element of $B$. Put $q(B):=\sum_{\alpha\in\Zer
f} O(B,\alpha)$. Then in the case when $\gamma$ grows from $B$ we
have
\begin{equation}\label{Eq3}
    \ord f(x,\gamma(x)) =
    \sum_{\alpha\in\Zer f} O(\gamma,\alpha) =
    \sum_{\alpha\in\Zer f} O(B,\alpha) = q(B).
\end{equation}

\medskip

\noindent The next lemma is a reformulation of Lemma 3.3 in
\cite{Kuo-Lu} (see also Lemma 2.2 in \cite{GP1}).

\begin{Lemma}\label{L1}\textbf{(The Kuo-Lu Lemma)} \\
(i)\quad For every $\gamma\in\Zer f_y'$ there exists $B\in T(f)$ of
          finite height such that $\gamma$ grows from $B$.\\
(ii)\quad For given $B\in T(f)$ of finite height the number of
         Puiseux roots of $f_y'$ which grow from $B$ counted with
        multiplicities is equal to $t(B)-1$ where $t(B)$ is the number
       of children of $B$ in $T(f)$.
\end{Lemma}

\noindent \textbf{Proof.}

\noindent Recall that $\Zer f=\{\alpha_1,\ldots,\alpha_d\}$ and let
$\Zer f'_y:=\{\gamma_1,\ldots,\gamma_{d-1}\}$. According to Lemma
3.3 in \cite{Kuo-Lu}, for given $\alpha\in \Zer f$ and a positive
rational number $r$,
\begin{equation}
\label{KL} \sharp \{j\;:\;O(\alpha,\gamma_j)=r\}= \sharp
\{k\;:\;O(\alpha,\alpha_k)=r\}.
\end{equation}
\noindent For the proof of the first statement fix $\gamma\in\Zer
f_y'$ and let $\alpha\in\Zer f$ be such that
$O(\alpha,\gamma)=\hbox{\rm max}_kO(\alpha_k,\gamma)$. Then from
(\ref{KL}) there exists $\alpha'$ such that
$O(\alpha,\gamma)=O(\alpha,\alpha')$ and consequently $\gamma$ grows
from the pseudo-ball $B(\alpha,O(\alpha,\alpha'))$.

\medskip
\noindent For the proof of the second statement suppose that
$B_1,\ldots,B_k$ are the children of $B_0$. By (\ref{KL}) for every
$B\in T(f)$ we have $\sharp (B\cap \Zer f)=\sharp
\{j\;:\;\gamma_j\in B\}+1$. Hence
\begin{eqnarray*}
\sharp \{j\;:\;\gamma_j\in\Zer f'_y \;\hbox{\rm grows from }B_0\}&=&
\sharp \{j\;:\;\gamma_j\in B_0\} - \sum_{i=1}^k\sharp
\{j\;:\;\gamma_j\in B_i\}\\&=& (\sharp (B_0\cap \Zer
f)-1)-\sum_{i=1}^k (\sharp (B_i\cap \Zer f)-1)\\&=&k-1,
\end{eqnarray*}
\noindent since $B_0\cap \Zer f=\displaystyle \bigsqcup_{i=1}^k
(B_i\cap \Zer f)$.

\medskip

\noindent The Kuo-Lu lemma together with (\ref{Eq3}) gives a
complete information on polar invariants and their multiplicities in
terms of $T(f)$. More precisely if $\tilde{\cal B}:=\{B\in
T(f)\;:\,h(B)<+\infty\}$ then the jacobian Newton polygon of
$f(x,y)=0$ equals

$$\sum_{B\in \tilde{\cal B}}
\left\{\Teisssr{(t(B)-1)q(B)} {t(B)-1}{18}{9}\right \}.$$

\begin{Example}
\label{Ex1}
Let $f(x,y)=(y^3-x^5)\prod_{i=1}^3(y-a_ix^2)$, $a_i\neq
a_j$ for $i\neq j$. The Puiseux roots of $f$ are: $\alpha_i=a_ix^2$
for $i=1,2,3$ and $\alpha_i=\epsilon^{i}x^{5/3}$ for $i=4,5,6$ and
$\epsilon=e^{2\pi \mathbf{i}/3}$.

\noindent Let us draw the Kuo-Lu tree of $f$.  Following
\cite{Kuo-Lu} we draw pseudo-balls of finite height as horizontal
bars and we do not draw pseudo-balls of infinite height. Three
Puiseux roots of $f_y'$ grow from a bar $B_1$ of height $5/3$ and
two Puiseux roots of a partial derivative grow from a bar $B_2$ of
height $2$. Since $q(B_1)=6\cdot(5/3)=10$ and
$q(B_2)=3\cdot2+3\cdot(5/3)=11$ we have $Q(f)=\langle\, 10:3,11:2
\,\rangle$.

\noindent Now take an irreducible series
$g(x,y)=(y^3-x^5)^2-9x^{11}$. Its Puiseux roots are
$\alpha_i(x)=\epsilon^{10i}x^{5/3}+\epsilon^{13i}x^{13/6}+\dots$ for
$i=1,\dots,6$ where $\epsilon$ is $6$-th primitive root of unity and
dots mean terms of
 higher degrees.

\noindent The tree $T(g)$ has one bar $B_1$ of height $5/3$ and
three bars $B_i$ of
 height~$13/6$.
Two Puiseux roots of $f_y'$ grow from $B_1$ and
$q(B_1)=6\cdot(5/3)=10$. From every $B_i$ ($i=2,3,4$) grows exactly
one Puiseux root of $f_y'$ and $q(B_i)=2\cdot(13/6)+4\cdot(5/3)=11$
for $i=2,3,4$. Hence $Q(g)=\langle\, 10:2,11:3 \,\rangle$.
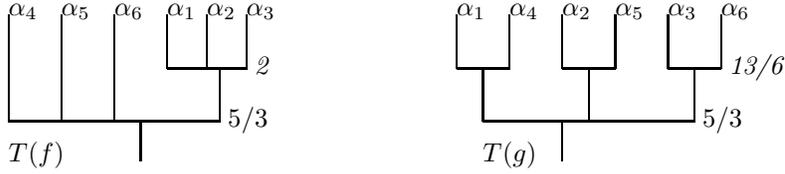
$$\mbox{%
\begin{picture}(100,60)(0,0)
\put(0,0){$T(f)$} \put(50,0){\line(0,1){15}} {\thicklines
\put(0,15){\line(1,0){80}}} \put(83,13){$5/3$}
\put(0,15){\line(0,1){40}} \put(20,15){\line(0,1){40}}
\put(40,15){\line(0,1){40}} \put(80,15){\line(0,1){20}} {\thicklines
\put(60,35){\line(1,0){30}}} \put(93,33){2}
\put(60,35){\line(0,1){20}} \put(75,35){\line(0,1){20}}
\put(90,35){\line(0,1){20}} \put(0,55){$\alpha_4$}
\put(20,55){$\alpha_5$} \put(40,55){$\alpha_6$}
\put(60,55){$\alpha_1$} \put(75,55){$\alpha_2$}
\put(90,55){$\alpha_3$}
\end{picture}
\hspace{1in}
\begin{picture}(100,60)(0,0)
\put(0,0){$T(g)$} \put(30,0){\line(0,1){15}} {\thicklines
\put(0,15){\line(1,0){80}}} \put(83,13){$5/3$}
\put(0,15){\line(0,1){20}} \put(40,15){\line(0,1){20}}
\put(80,15){\line(0,1){20}} {\thicklines
\put(-10,35){\line(1,0){20}}}
\put(-10,35){\line(0,1){20}}\put(-10,55){$\alpha_1$}
\put(10,35){\line(0,1){20}}\put(10,55){$\alpha_4$} {\thicklines
\put(30,35){\line(1,0){20}}}
\put(30,35){\line(0,1){20}}\put(30,55){$\alpha_2$}
\put(50,35){\line(0,1){20}}\put(50,55){$\alpha_5$} {\thicklines
\put(70,35){\line(1,0){20}}}
\put(70,35){\line(0,1){20}}\put(70,55){$\alpha_3$}
\put(90,35){\line(0,1){20}}\put(90,55){$\alpha_6$} \put(93,33){13/6}
\end{picture}}
$$
\end{Example}

\begin{Corollary}\label{Cor1}
Let $f=f(x,y)\in\Cn\{x,y\}$ be a power series such that the vertical
axis is transverse to~$\{f=0\}$ and let $w=w(x,y)\in\Cn\{x\}[y]$ be
the Weierstrass polynomial of $f$. Then ${\cal N}_J(f)={\cal
N}_J(w)$.
\end{Corollary}

\noindent \textbf{Proof.} Write $f(x,y)$ in the form~(\ref{Eq1}).
Then $w(x,y)=\prod_{i=1}^d(y-\alpha_i(x))$ and clearly $\Zer f=\Zer
\,w$. Since polar invariants and their multiplicities depend only on
Kuo-Lu tree we have ${\cal N}_J(f)={\cal N}_J(w)$.

\section{Similarity lemma}
Let $f$ be an irreducible power series in two variables of order
greater than 1. Then after an analytic change of coordinates we may
assume that

\begin{equation}\label{Eq4}
    f(x,y)= (y^n-x^m)^a + \sum_{ni+mj>anm} f_{i,j}x^iy^j
\end{equation}

\noindent where $1<n<m$ are coprime integers. By Merle's formula
(see \cite{Merle},\cite{GP1}) the smallest polar invariant of~$f$
is~$am$ with multiplicity~$n-1$. We will show that the converse is
also true.

\begin{Lemma}\label{L2}
Let $f\in\Cn\{x,y\}$ be a power series such that $\ord f=an$ and the
smallest polar invariant of~$f$ is~$am$ with multiplicity~$n-1$. If
$1<n<m$ are coprime integers, then after an analytic change of
coordinates
$$  f(x,y)= (y^n-x^m)^a + \sum_{ni+mj>anm} f_{i,j}x^iy^j .
$$
\end{Lemma}

\noindent \textbf{Proof.} Choose a system of coordinates such that
the $y$-axis is transverse to~$\{f=0\}$ and write a Puiseux
factorization of $f$
$$ f(x,y)=\mathit{unit}\prod_{i=1}^{an}(y-\alpha_i(x)).
$$
Let $B$ be the root of $T(f)$. Since contact orders
$O(\alpha_i,\alpha_j)$ are greater than or equal to $h(B)$ all
Puiseux roots of $f$ have a form
$\alpha_i(x)=\lambda(x)+c_ix^{h(B)}+\dots$ where
$\lambda(x)=a_1x^{\delta_1} + a_2x^{\delta_2} + \cdots +
a_kx^{\delta_k}$ ($1\leq \delta_1<\delta_2<\dots<\delta_k$) is a
finite sum of terms of degrees smaller than $h(B)$ and at least one
$c_i$ is non zero. We will show that all exponents in $\lambda(x)$
are integers. Suppose to the contrary that $\delta_j$ are integers
for $1\leq j<s$ and $\delta_s=p/q$, for $p$, $q$ coprime integers.
Then by [Walk, p.~107 
Theorem~4.1], if $\omega$ is a $q$-th primitive root of unity, a
series $\bar\alpha$ of the form
$$ \bar\alpha(x)= a_1x^{\delta_1} + \cdots + a_{s-1}x^{\delta_{s-1}} +
   \omega^pa_sx^{p/q}+\dots,
$$
which is a conjugate of
$$ \alpha(x)= a_1x^{\delta_1} + \cdots + a_{s-1}x^{\delta_{s-1}}
+ a_sx^{p/q}+\dots
$$
is a Puiseux root of $f$ and we get a contradiction because
$O(\alpha,\bar\alpha)=p/q<h(B)$.

\medskip

\noindent We checked that $\lambda(x)$ is a polynomial. After an
analytic substitution $y:=y-\lambda(x)$ we may assume that
$$ \alpha_i(x) = c_ix^{h(B)}+\dots \qquad\mbox{for $i=1,\dots,an.$}
$$

\noindent Computing the smallest polar invariant of $f$ by
formula~(\ref{Eq3}) we get
$$ am = q(B) = \sum_{\alpha\in \Zer f} O(B,\alpha) =
          \sum_{\alpha\in \Zer f} h(B) = an\,h(B) .
$$

\noindent Hence $h(B)=m/n$.

\medskip

\noindent Fix a weight $w$ such that $w(x)=n$, $w(y)=m$ and let
$\ini f$ denotes the initial quasi-homogeneous part of $f$ with
respect to~$w$. Since all factors $y-c_ix^{m/n}$ are quasi
homogeneous we have
\begin{equation}\label{Eq5}
 \ini f(x,y) = \mathit{const}\prod_{i=1}^{an}(y-c_ix^{m/n}).
\end{equation}

\noindent On the other hand because $\ini f(x,y)$ is a
quasi-homogeneous polynomial it can be written as
$$ \mathit{const}\prod_{j=1}^s (y^n-C_jx^m)^{k_j} \;\;\;
$$
where $C_i\neq C_j$ for $i\neq j$ and $k_1+\dots+k_s=a$. We want to
show that $s=1$. Suppose to the contrary that $s\geq 2$. Assume for
simplicity that $C_1\neq 0$. Then $(y^n-C_1x^m)(y^n-C_2x^m)$ equals
\begin{eqnarray*}
      \;\;\prod_{j=1}^n (y-\omega^j\sqrt[n]{C_1}x^{m/n})(y-\omega^j
      \sqrt[n]{C_2}x^{m/n})
    & \quad\mbox{if $C_2\neq0$} \\
      y^n\prod_{j=1}^n (y-\omega^j\sqrt[n]{C_1}x^{m/n})
    & \quad\mbox{if $C_2=0$}
\end{eqnarray*}

\noindent where $\omega$ is an $n$-th primitive root of unity. Hence
(\ref{Eq5}) has at least $n+1$ different factors. Because different
$c_i$ yield different edges of $T(f)$ growing from the root of
$T(f)$, the multiplicity of the smallest polar invariant is at
least~$n$ and we arrive at contradiction.

\medskip

\noindent We checked that $\ini f(x,y)=\mathit{const}(y^n-C_1x^m)^a$
and finally a substitution $x:=C_1^{1/m}x$ gives $C_1=1$.

\section{Reduction technique}
\label{reduction technique}

\noindent Take a distinguished Weierstrass polynomial
$f(x,y)\in\Cn\{x\}[y]$ of the form~(\ref{Eq4}) not necessarily
irreducible. Let $F(x,y)=f(x^n,y)$ and let $\omega$ be an $n$-th
primitive root of unity. By equality
$$ (y^n-x^{nm})^a=(y-x^m)^a(y-\omega x^m)^a\cdots(y-\omega^{n-1} x^m)^a
$$
and Hensel's lemma we get a factorization $F=f_0f_1\cdots f_{n-1}$
where $f_i$ are Weierstrass polynomials with quasi-homogeneous
initial forms $(y-\omega^i x^m)^a$ for $i=0,\dots,n-1$. Denote by
$\tilde f$ the factor $f_0$. This transformation has nice
properties.

\begin{Lemma} \ \label{L3}
\begin{description}
\item[(i)]
  $f$ is irreducible if and only if $\tilde f$ is irreducible.

\item[(ii)]
  Let $Q(f)=\langle\, q_0:m_0,q_1:m_1,\dots,q_r:m_r\, \rangle$
  be the system of polar invariants of $f$ with $q_{i-1}<q_i$ for $i=1,\dots,r$.
  Then
  $Q(\tilde f)=\langle\, q_1':m_1',\dots,q_r':m_r'\, \rangle$ where
  \begin{eqnarray*}
        m_i' &=& m_i/(m_0+1) \\
        q_i' &=& q_i(m_0+1)-q_0m_0
  \end{eqnarray*}
  for $i=1,\dots,r$.

\item[(iii)] If f is irreducible with Puiseux characteristic
  $(\beta_0,\ldots,\beta_r)$ then
  the Puiseux characteristic of $\tilde f$ is
  $\left(\frac{\beta_0}{n},\beta_2,\ldots,\beta_r\right)$.
\end{description}
\end{Lemma}

\noindent \textbf{Proof.} Assume that $f$ is irreducible. Then all
Puiseux roots $\alpha_i$ of $f$ are fractional power series and the
exponents of those series have a least common denominator $na$.
Since Puiseux roots of $F(x,y)=f(x^n,y)$ are $\alpha_i(x^n)$ the
exponents of these series have a least common denominator $a$. Thus
the order of irreducible factors of $F$ is exactly $a$ and we see
that $F=f_0f_1\cdots f_{n-1}$ is a decomposition into irreducible
factors.

\medskip

\noindent If $f=h_1h_2$ is a product then clearly $\tilde f=\tilde
h_1 \tilde h_2$. Thus we ended proof of (i).

\vspace{2ex} \noindent To prove~(ii) we will compute the system of
polar invariants~$Q(F)$ in two ways. Write Puiseux factorizations
of~$f$ and~$f_y'$:
$$    f(x,y)=  \prod_{i=1}^{an}(y-\alpha_i(x)), $$
$$ f_y'(x,y)=an\prod_{j=1}^{an-1}(y-\gamma_j(x)). $$

\noindent The Puiseux roots of $F$ are
$\bar\alpha_i(x)=\alpha_i(x^n)$ for $i=1,\dots, an$ and by equality
$\frac{\partial F}{\partial y}(x,y)=\frac{\partial f}{\partial
y}(x^n,y)$, the Puiseux roots of $F_y'$ are
$\bar\gamma_i(x)=\gamma_i(x^n)$ for $i=1,\dots, an-1$. We have
\begin{equation}\label{Eq7}
\ord F(x,\bar\gamma_i(x))=\ord f(x^n,\gamma_i(x^n))=n\,\ord
f(x,\gamma_i(x))
\end{equation}
for $i=1,\dots,an-1$. Hence
\begin{equation}\label{EqQ1}
Q(F)=\langle\, nq_0:m_0, \dots , nq_r:m_r \,\rangle.
\end{equation}

\medskip

\noindent In the second part of computation we will show that the
Kuo-Lu tree-model of $F=f_0\cdots f_{n-1}$ has a special structure.
It separates above the root to $n$ ~sub-trees and each of them is
isomorphic to $T(\tilde f)$. This together with Kuo-Lu Lemma will
allow to express the system of polar invariants of $F$ by the system
of polar invariants of $\tilde f$.

\medskip

\noindent By \cite{Walker} if $\bar\alpha\in\Zer F$ is a Puiseux
root of $f_i$ then the initial term of $\bar\alpha$ is
$\omega^ix^m$. Hence for every $\bar\alpha$, $\bar\alpha'\in\Zer F$
\begin{eqnarray*}
       O(\bar\alpha,\bar\alpha')=m &
       \mbox{if $\bar\alpha\in\Zer f_i$, $\bar\alpha'\in\Zer f_j$, $i\neq j$} \\
       O(\bar\alpha,\bar\alpha')>m &
       \mbox{if $\bar\alpha$, $\bar\alpha'\in\Zer f_i$}
\end{eqnarray*}

\noindent It follows that the Kuo-Lu tree-model $T(F)$ has a root
$B_0$ of height $m$ and above the root it separates to $n$
sub-tree-models $T_i$ ($i=0,\dots,n-1$); the leaves of $i$-th
sub-tree-model are roots of $f_i$. Moreover if $\phi$ grows
from~$B\in T_i$ (in short $\phi$ grows from~$T_i$) then
$\phi=\omega^ix^m+\dots$.

\medskip

\noindent We shall establish a one-to-one correspondence between the
Puiseux roots of~$F_y'$ which grow from $T_0$ and those which grow
from $T_i$ ($i=1,\dots,n-1$). To do this we need some properties of
action of complex roots of unity on fractional power series.

\medskip

\noindent Let $D$ be an integer such that every $\phi\in\Zer
f\cup\Zer f_y'$ can be written as fractional power series with
common denominator $D$.
$$ \phi(x)=a_1x^{n_1/D}+a_2x^{n_2/D}+\dots, \qquad
   1\leq D<n_1<n_2<\dots
$$

\noindent Clearly $D$ is a multiple of $n$ because $m/n$ is the
smallest exponent of $\phi\in\Zer f$. For every $\theta\in\Cn$ such
that $\theta^D=1$ we define a conjugate $\theta(\phi)$ by formula
$$ \theta(\phi)(x)=a_1\theta^{n_1}x^{n_1/D}+a_2\theta^{n_2}x^{n_2/D}+\dots
$$

\noindent It is shown in \cite{Kuo-Pa} (page~293) that for every
$D$-th root of unity $\theta$ the action of~$\theta$ permutes the
sets $\Zer f$ and $\Zer f_y'$ and preserves a contact order, that is
$O(\phi_1,\phi_2)=O(\theta(\phi_1),\theta(\phi_2))$.

\medskip

\noindent \textbf{Claim 1.} $\ord f(x,\gamma(x))=\ord
f(x,\theta(\gamma)(x))$ for every $\gamma\in\Zer f_y'$,
$\theta^D=1$.

\vspace{2ex} \noindent The proof is straightforward. From properties
of conjugation mentioned above we get $\ord f(x,\gamma(x)) =
 \sum_i O(\gamma,\alpha_i) =
 \sum_i O(\theta(\gamma),\theta(\alpha_i)) =
 \sum_i O(\theta(\gamma),\alpha_i) = \ord f(x,\theta(\gamma)(x))
$.

\vspace{2ex} \noindent Recall that $\omega$ is a primitive $n$-th
root of unity.

\medskip

\noindent \textbf{Claim 2.} There exist $\theta$ ($\theta^D=1$) such
that for every $\gamma\in\Zer f_y'$ of the form
$\gamma(x)=x^{m/n}+\dots$ we have
$$ \theta^i(\gamma)(x)=\omega^ix^{m/n}+\dots \quad\mbox{for $i=1,\dots,n-1$} .
$$

\noindent \textit{Proof of Claim 2.} Since $m$ and $n$ are coprime
there is $pm+qn=1$ for some integers~$p$, $q$. Let $\theta$ be the
complex number such that $\theta^{D/n}=\omega^p$. Clearly
$\theta^D=1$. We have
$(\theta^i)^{m(D/n)}=\omega^{imp}=\omega^{imp}\omega^{inq}=
 \omega^{imp+inq}=\omega^i$.
Writing $\gamma$ as a power series with common denominator $D$
$$ \gamma(x)=x^{\frac{m(D/n)}{D}}+\dots
$$

\noindent we see that
$\theta^i(\gamma)(x)=(\theta^i)^{m(D/n)}x^{\frac{m(D/n)}{D}}+\dots=\omega^ix^{m/n}+\dots$

\vspace{2ex} \noindent Observe that for every $\gamma\in\Zer f_y'$ a
series $\gamma(x^n)$ grows from a sub-tree-model~$T_i$ if and only
if $\gamma$ has a form $\gamma(x)=\omega^ix^{m/n}+\dots$ for
$i=0,\dots,n-1$. Thus by Claim~2 the action of $\theta^i$ yields a
one-to-one correspondence between Puiseux roots of $F_y'$ which grow
from~$T_0$ and Puiseux roots of $F_y'$ which grow from~$T_i$.
By~(\ref{Eq7}) and~Claim~1 the corresponding roots give the same
polar invariants.

\medskip

\noindent Take $\bar\gamma\in\Zer F_y'$.  If $\bar\gamma$ grows from
a root $B_0$ then $O(\bar\gamma,\bar\alpha)=m$ for all
$\bar\alpha\in\Zer F$. Hence $\ord
F(x,\bar\gamma(x))=\sum_{\bar\alpha\in\Zer
f_i}O(\bar\gamma,\bar\alpha)=anm$. Because from $B_0$ grow $n$
sub-tree-models there are $n-1$ Puiseux roots of $F_y'$ which give
the smallest polar invariant $anm$.

\medskip

\noindent Now take $\bar\gamma(x)=x^m+\dots$ which grows from $T_0$.
Then for $i=1,\dots,n-1$ we have $\ord f_i(x,\bar\gamma(x)) =
\sum_{\bar\alpha\in\Zer f_i} O(\bar\gamma,\bar\alpha) = am$. By
equality $\ord F(x,\bar\gamma(x)) = \sum_{i=0}^{n-1} \ord
f_i(x,\bar\gamma(x))$ we get
\begin{equation}\label{Eq8}
\ord F(x,\bar\gamma(x))=\ord \tilde f(x,\bar\gamma(x))+(n-1)am .
\end{equation}

\noindent The sub-tree-model $T_0$ is equal to $T(\tilde f)$. By the
Kuo-Lu lemma there is one-to-one correspondence between Puiseux
roots of $\tilde f_y'$ and Puiseux roots of $F_y'$  which grow from
$T_0$. Thus if $Q(\tilde f)=\langle\, q_1':m_1' \dots ,q_s':m_s'
\,\rangle$ then from (\ref{Eq8})
\begin{equation}\label{EqQ2}
Q(F)=\langle\, amn:n-1,q_1'+(n-1)am:nm_1' \dots ,q_s'+(n-1)am:nm_s'
\,\rangle.
\end{equation}

\noindent Comparing~(\ref{EqQ1}) and~(\ref{EqQ2}) we get (ii).

\medskip

\noindent To prove (iii) recall that the Puiseux roots of $\tilde f$
are exactly the Puiseux roots $\bar{\alpha}_i$ of $F(x,y)=0$ such
that $\hbox{\rm in}(\bar{\alpha}_i)=x^m$. Put $\hbox{\rm Zer}(\tilde
f)=\{\bar{\alpha}_{j_1},\ldots,\bar{\alpha}_{j_a}\}$ then
$$ \{\hbox{\rm ord}(\bar{\alpha}_{j_k}-\bar{\alpha}_{j_l})\}_{k\neq l}=
\{n.\hbox{\rm ord}({\alpha}_{j_k}-{\alpha}_{j_l})\}_{k\neq l}=
\left\{n\frac{\beta_2}{\beta_0},\ldots,n\frac{\beta_g}{\beta_0}\right\}=
\left\{\frac{\beta_2}{a},\ldots,\frac{\beta_g}{a}\right\}$$

\noindent so $\hbox{\rm char}(\tilde f)=(a,\beta_2,\ldots,\beta_g)$.

\medskip

\vspace{2ex}

\noindent \textbf{Proof of Theorem \ref{Th1}.} We apply induction
with respect to $\ord f$. If $\ord f=1$ then $f$ has an empty system
of polar invariants. Hence also $g_y'$ has no Puiseux root and
consequently $\ord g=1$. In this case both $f$ and $g$ are
irreducible.

\noindent Now assume that $\ord f>1$ and that Theorem is true for
smaller orders. We may assume that $f$ has form (\ref{Eq4}). From
${\cal N}_J(f)={\cal N}_J(g)$ it follows that $g$ satisfies the
assumptions of Lemma~\ref{L2} so we may also assume that~$g$ has
form (\ref{Eq4}). By Corollary~\ref{Cor1} we may replace $f$ and $g$
by their distinguished Weierstrass polynomials. Take $\tilde f$ and
$\tilde g$. By Lemma~\ref{L3} ${\cal N}_J(\tilde f)={\cal
N}_J(\tilde g)$ and $\tilde f$ is irreducible. Hence by inductive
hypothesis $\tilde g$ is irreducible so also $g$ is.

\begin{Remark}

\noindent If we put $f=y^2(y-x^2)^2+x^{11}$ and
$g=y^3(y-x^2)+x^{11}$, then ${\cal N}_J(f)={\cal
N}_J(g)=\{\Teis{8}{1}\}+\{\Teis{22}{2}\}$ but $f=0$ has two branches
while $g=0$ has four. So we cannot generalize Theorem 1 to
multi-branched curves.
\end{Remark}

\begin{Remark}
In \cite{Lenarcik} the author gives an example of two polynomials
$f(x,y)=(y-x^2)(y^4-x^{12})$ and $g(x,y)=(y-x^2)^2y^3+x^{14}$ such
that the curves $f=0$ and $g=0$ are unitangent and have the same
jacobian Newton polygon, but $f=0$ is non-degenerate with 3
irreducible components and $g=0$ is degenerate with 5 irreducible
components.
\end{Remark}

\begin{Remark}
In Example \ref{Ex1} we take two curves $f=0$ and $g=0$, where $f=0$
has 4 branches and $g=0$ is irreducible. These curves have the same
set of polar invariants, so in Theorem \ref{Th1} we cannot replace
the jacobian Newton polygon by the collection of inclinations of its
edges.
\medskip

\noindent Merle in \cite{Merle} proves that the set of polar
invariants and the order of a branch determine its equisingularity
class. Nevertheless this data is the same for the curves $f=0$ and
$g=0$ of Example \ref{Ex1}, so the knowledge of this data is not
enough to decide if a curve is irreducible or not.
\end{Remark}

\medskip \noindent Theorem \ref{Th1} allows us to decide if a curve
is irreducible using its jacobian Newton polygon. In the sequel,
following this approach, we propose new criteria of irreducibility
of complex plane curves.

\section{Characterization of jacobian Newton polygons of branches
using the reduction operation}

\medskip

\noindent Note that if $\tilde f$ is the reduction of $f$, then by
(ii) of Lemma \ref{L3}
$${\cal N}_J(\tilde f):=\sum_{j=1}^{r}
\left\{\Teisssr{m'_jq'_j}{m'_j}{6}{3}\right \}.$$

\medskip

\noindent We will call the transition from ${\cal N}_J(f)$ to ${\cal
N}_J(\tilde f)$ a reduction operation. In this section we extend
this operation to all rational Newton polygons and using the results
of section \ref{reduction technique} we will characterize the
special convenient Newton polygons which are jacobian Newton
polygons of an irreducible plane curve.

\vspace{0.7cm} \noindent {\bf The reduction operation over rational
Newton polygons:}

\medskip

\noindent Let ${\cal
N}=\sum_{i=1}^r\left\{\Teisssr{L_i}{M_i}{4}{2}\right \}$ be a
convenient rational Newton polygon with $r\geq 2$ and
$\frac{L_1}{M_1}<\cdots < \frac{L_r}{M_r}$. The {\em reduction} of
${\cal N}$ is by definition
\begin{equation}
\label{reduction polygon} {\cal R}({\cal
N}):=\sum_{i=1}^{r-1}\left\{\Teisssr{L'_i}{M'_i}{4}{2}\right\}
\end{equation}
\noindent where

\medskip

\begin{equation}
\label{def1} L'_i:=L_{i+1}-\frac{L_1}{1+M_1}M_{i+1},
\end{equation}

\medskip

\noindent and

\begin{equation}
\label{def2}
 M'_i:=\frac{M_{i+1}}{1+M_1}.
\end{equation}
\bigskip

\noindent Since $\frac{L_{i+1}}{M_{i+1}}<\frac{L_{i+2}}{M_{i+2}}$,
${\cal R}({\cal N})$ is again a convenient rational Newton polygon
written in its canonical form because
$\frac{L'_i}{M'_i}<\frac{L'_{i+1}}{M'_{i+1}}$.

\medskip

\noindent The reduction operation transforms a Newton polygon of
$r>1$ compact faces to a Newton polygon of $r-1$ compact faces.

\medskip

\noindent We denote by ${\cal R}^i({\cal N})$ the Newton polygon
obtained after applying $i$-times the reduction operation to ${\cal
N}$. By convention we put ${\cal R}^0({\cal N})={\cal N}$.

\medskip

\noindent If ${\cal
N}=\sum_{k=1}^r\left\{\Teisssr{L_k}{M_k}{4}{2}\right \}$ then we
denote ${\cal R}^i({\cal N})$ by
$\sum_{k=1}^{r-i}\left\{\Teisssr{L^{(i)}_k}{M^{(i)}_k}{6}{3}\right
\}$.

\medskip

\begin{Properties}
\label{properties} $\;$
\begin{enumerate}
\item
${\cal R}({\cal N}_J(f))={\cal N}_J(\tilde f)$.
\item If ${\cal N}=\sum_{i=1}^r\left\{\Teisssr{L_i}{M_i}{4}{2}\right
\}$ and ${\cal
N}^*=\sum_{i=1}^r\left\{\Teisssr{L^*_i}{M^*_i}{4}{2}\right \}$
satisfy ${\cal R}({\cal N})={\cal R}({\cal N}^*)$ and
$\left\{\Teisssr{L_1}{M_1}{4}{2}\right
\}=\left\{\Teisssr{L^*_1}{M^*_1}{4}{2}\right \}$ then ${\cal
N}={\cal N}^*$.
\end{enumerate}
\end{Properties}

\noindent \textbf{Proof.} Using Lemma \ref{L3} we have (1). The
property (2) follows from (\ref{def1}) and (\ref{def2}) since
$M_{i+1}=M^*_{i+1}$ and  $L_{i+1}=L^*_{i+1}$.

\bigskip

\begin{Theorem}
\label{second main theorem} Let ${\cal
N}=\sum_{k=1}^r\left\{\Teisssr{L_k}{M_k}{4}{2}\right \}$ be a
special convenient integral Newton polygon satisfying

\begin{enumerate}
\item $1+\hbox{\rm ht}({\cal N})<\frac{L_1}{M_1}$,

\item ${\cal R}^i({\cal N})$ is an integral Newton polygon and
$\frac{L_1^{(i)}}{M_1^{(i)}} \in \mathbf {N}$ for  $0\leq i\leq
r-1,$


\item $(1+M_1^{(i)})\hbox{\rm g.c.d.}\left(\frac{L_1^{(i)}}{M_1^{(i)}},
1+\hbox{\rm ht}({\cal R}^i({\cal N}))\right)=1+\hbox{\rm ht}({\cal
R}^i({\cal N}))$ for $0\leq i\leq r-1$,
\end{enumerate}

\noindent then there exists $f\in\mathbf C\{x,y\}$ irreducible such
that ${\cal N}={\cal N}_{J}(f)$ and the Puiseux characteristic of
$f(x,y)=0$ is $\left(1+\hbox{\rm ht}({\cal
N}),\frac{L_1^{(0)}}{M_1^{(0)}}, \frac{L_1^{(1)}}{M_1^{(1)}},\ldots,
\frac{L_1^{(r-1)}}{M_1^{(r-1)}}\right)$.

\end{Theorem}

\noindent \textbf{Proof.-}

\noindent We use induction on $r$. According to Merle (\cite{Merle})
the Newton polygon ${\cal N}=\left\{\Teisssr{L_1}{M_1}{4}{2}\right
\}$ is the jacobian Newton polygon of a branch $f(x,y)=0$ if and
only if $\hbox{\rm g.c.d.}(\frac{L_1}{M_1},1+M_1)=1$ and
$1+M_1<\frac{L_1}{M_1}$, moreover in this case the Puiseux
characteristic of $f(x,y)=0$ is $(1+M_1,\frac{L_1}{M_1})$.

\medskip

\noindent We suppose now that Theorem \ref{second main theorem} is
true for any special convenient integral Newton polygon with $r-1$
compact edges and let ${\cal
N}=\sum_{k=1}^r\left\{\Teisssr{L_k}{M_k}{4}{2}\right \}$ be a
special convenient integral Newton polygon satisfying the three
conditions of Theorem \ref{second main theorem}.  Then

\medskip

\noindent \noindent \textbf{Claim 1.} The reduction ${\cal R}({\cal
N})$ of ${\cal N}$ also satisfies the hypothesis of Theorem
\ref{second main theorem}.

\medskip
\noindent The second and the third conditions are clearly satisfied.
By hypothesis $\frac{L_1}{M_1}>1+\hbox{\rm ht}({\cal N})$ so
\begin{eqnarray}
\label{T21} \frac{L_2}{M_2}>\frac{L_1}{M_1}>1+\hbox{\rm ht}({\cal
N})
>\frac{1+\hbox{\rm ht}({\cal N})}{1+M_1}=1+\hbox{\rm ht}({\cal R}({\cal N})).
\end{eqnarray}

\noindent Moreover
\begin{equation}
\label{T22} \frac{L'_1}{M'_1}>\frac{L_2}{M_2}\end{equation}

\noindent since $M_1L_2-L_1M_2>0$, then (\ref{T21}) and (\ref{T22})
give $\frac{L'_1}{M'_1}>1+\hbox{\rm ht}({\cal R}({\cal N})).$

\bigskip

\noindent So there is an irreducible plane curve $f'$ such that
${\cal N}_J(f')={\cal R}({\cal N})$.\\ Put
$(\beta^{'}_0,\beta^{'}_1,\ldots,\beta^{'}_{r-1})$ the Puiseux
characteristic of $f'=0$. We define $n$ and $m$ from ${\cal N}$ in
the following way:

\begin{enumerate}
\item $n:=1+M_1$,
\item $m:=\frac{\frac{L_1}{M_1}}{\hbox{\rm g.c.d}(1+\hbox{\rm
ht}({\cal N}),\frac{L_1}{M_1})}$.
\end{enumerate}

\medskip

\noindent \noindent \textbf{Claim 2.} The numbers $n$ and $m$ are
coprime integers and $m>n$.

\medskip

\noindent Put $a:=\hbox{\rm g.c.d}(1+\hbox{\rm ht}({\cal
N}),\frac{L_1}{M_1})$. By definition $\frac{L_1}{M_1}=ma$ and for
$i=0$ the third condition of the theorem gives $a=\hbox{\rm
g.c.d.}(na,ma)$, then $n$ and $m$ are coprime. Moreover $m>n$ since
$\frac{L_1}{M_1}>1+\hbox{\rm ht}({\cal N}).$

\medskip

\noindent \textbf{Claim 3.} The numbers
$(n\beta^{'}_0,m\beta^{'}_0,\beta^{'}_1,\ldots,\beta^{'}_{r-1})$
form an increasing sequence of coprime integers.

\medskip

\noindent Observe that

\begin{eqnarray*}
\hbox{\rm
g.c.d.}(n\beta^{'}_0,m\beta^{'}_0,\beta^{'}_1,\ldots,\beta^{'}_{r-1})&=&
\hbox{\rm g.c.d.}(\hbox{\rm
g.c.d.}(n\beta^{'}_0,m\beta^{'}_0),\beta^{'}_1,\ldots,\beta^{'}_{r-1}))\\
&=&\hbox{\rm
g.c.d.}(\beta^{'}_0,\beta^{'}_1,\ldots,\beta^{'}_{r-1})=1.
\end{eqnarray*}

\noindent In order to prove that
$(n\beta^{'}_0,m\beta^{'}_0,\beta^{'}_1,\ldots,\beta^{'}_{r-1})$ is
an increasing sequence, it is enough to prove that
$m\beta^{'}_0<\beta^{'}_1$ that is equivalent to prove
$\frac{L_1}{M_1}<\frac{L'_1}{M'_1}$, which is true by the inequality
(\ref{T22}).

\bigskip

\noindent So there exists an irreducible plane curve $f(x,y)=0$
which Puiseux characteristic
$(n\beta^{'}_0,m\beta^{'}_0,\beta^{'}_1,\ldots,\beta^{'}_{r-1})$.

\bigskip

\noindent Consequently using Lemma \ref{L3} we obtain that
$\hbox{\rm char}(\tilde f)= \hbox{\rm char}(f')$ and then
$$\hbox{\rm
char}(f)=(n\beta'_0,m\beta'_0,\beta^{'}_1,\ldots,\beta^{'}_{r-1})=
\left(1+\hbox{\rm ht}({\cal N}),\frac{L_1^{(0)}}{M_1^{(0)}},
\frac{L_1^{(1)}}{M_1^{(1)}},\ldots,
\frac{L_1^{(r-1)}}{M_1^{(r-1)}}\right).
$$

\medskip

\noindent In order to finish we have to prove that ${\cal
N}_J(f)={\cal N}$. For that observe that $f'$ and $\tilde{f}$ are
irreducible with $\hbox{\rm char}(f')=\hbox{\rm char}(\tilde{f})$,
that is, $f'$ and $\tilde{f}$ are equisingular, so ${\cal
N}_J(f')={\cal N}_J(\tilde{f})$ and by (1) of Properties
\ref{properties} we can write ${\cal N}_J(f')={\cal R}({\cal
N}_J(f))={\cal R}({\cal N})$. By~(\ref{Me}) the first segment of the
jacobian Newton polygon ${\cal N}_J(f)$ has height $n-1=M_1$ and
inclination $m\beta'_0=L_1^{(0)}/M_1^{(0)}$ hence equality ${\cal
N}_J(f)={\cal N}$ follows from (2) of Properties \ref{properties}.

\medskip

\begin{Corollary}
\label{characterization} A special convenient integral Newton
polygon ${\cal N}=\sum_{k=1}^r\left\{\Teisssr{L_k}{M_k}{4}{2}\right
\}$ is the jacobian Newton polygon of a branch if and only if
verifies the next three conditions:

\begin{enumerate}

\item $1+\hbox{\rm ht}({\cal N})<\frac{L_1}{M_1}$,

\item ${\cal R}^i({\cal N})$ is an integral Newton polygon and
$\frac{L_1^{(i)}}{M_1^{(i)}}\in \mathbf N$ for all $0\leq i\leq
r-1,$

\item $(1+M_1^{(i)})\hbox{\rm g.c.d.}\left(\frac{L_1^{(i)}}{M_1^{(i)}},
1+\hbox{\rm ht}({\cal R}^i({\cal N}))\right)=1+\hbox{\rm ht}({\cal
R}^i({\cal N}))$ for $0\leq i\leq r-1$.
\end{enumerate}
\end{Corollary}

\noindent {\bf Proof.} It is enough to show that if ${\cal N}={\cal N}_J(f)$
for an irreducible distinguished Weierstrass polynomial $f(x,y)\in\Cn\{x\}[y]$
then ${\cal N}$ satisfies conditions 1--3.

\medskip

\noindent Let
$\bigl\langle\overline{\beta}_0,\overline{\beta}_1,\ldots,\overline{\beta}_g\bigr\rangle$
be the semigroup of $\{f=0\}$. Then by Merle formula~(\ref{Me})
$L_1=(n_1-1)\overline{\beta}_1$, $M_1=n_1-1$ where
$n_1=\overline{\beta}_0/\hbox{\rm
g.c.d.}(\overline{\beta}_0,\overline{\beta}_1)$ and $\hbox{\rm
ht}({\cal N})+1=\overline{\beta}_0$. Condition~1 reads as
$\overline{\beta}_0<\overline{\beta}_1$, condition~2 for $i=0$ is
satified because $\frac{L_1}{M_1}=\overline{\beta}_1$ and
condition~3 for $i=0$ is an equation $n_1\,\hbox{\rm
g.c.d.}(\overline{\beta}_1,\overline{\beta}_0)=\overline{\beta}_0$.
In order to check that conditions~2 and~3 are satisfied for $1\leq
i\leq g-1$ it is enough to observe that ${\cal R}({\cal N})={\cal
N}_J(\tilde f)$ and apply an induction because $\tilde f$ is again
an irreducible distinguished Weierstrass polynomial.

\medskip

\begin{Remark}
>From Corollary \ref{characterization} we have that if ${\cal N}={\cal N}_J(f)$,
with $f$ irreducible and $\hbox{\rm
char}(f)=(\beta_0,\ldots,\beta_g)$ then $\beta_0=1+\hbox{\rm
ht}({\cal N}_J(f))$ and $\beta_{k+1}=\frac{L_1^{(k)}}{M_1^{(k)}}$.
\end{Remark}

\section{Characterization of jacobian Newton polygons of branches
using the abrasion operation} \label{abrasion} \noindent The concept
of approximate root was introduced and studied in
\cite{Abhyankar-Moh}:

\begin{Proposition}
Let $A$ be an integral domain. If $f(y)\in \mathbf A[y]$ is monic of
degree $d$ and $p$ is invertible in $A$ and divides $d$, then there
is a unique monic polynomial $g(y) \in \mathbf A[y]$ such that the
degree of $f-g^p$ is less than $d-\frac{d}{p}$.
\end{Proposition}

\medskip

\noindent This allows us to define:

\begin{Definition}
The unique monic polynomial of the preceding proposition is named
the $p$-th approximate root of $f$.
\end{Definition}

\noindent Suppose now that $f\in \mathbf C\{x\}[y]$ is irreducible
of Puiseux characteristic $(\beta_0,\ldots,\beta_g)$. Put
$l_k:=\hbox{\rm g.c.d.}(\beta_0,\ldots,\beta_k)$. In particular
$l_k$ divides $\hbox{\rm deg}f=\beta_0$ for all $k\in
\{1,\ldots,g\}$. We note in what follows $f^{(k)}$ the $l_k$-th
approximate root of $f$ which we named {\em characteristic
approximate} roots of $f$.

\medskip

\noindent Next proposition is the main one in \cite{Abhyankar-Moh}
(see also \cite{GP2} and \cite{Popescu}):

\begin{Proposition}
Let $f$ be an irreducible curve with characteristic Puiseux
$(\beta_0,\ldots,\beta_g)$ and semigroup $\overline{\beta}_0\mathbf
N+\cdots+\overline{\beta}_g\mathbf N$. The approximate roots
$f^{(k)}$ of $f$, for $0\leq k \leq g$, have the following
properties:

\begin{enumerate}
\item The degree of $f^{(k)}$ is equal to $\frac{\beta_0}{l_k}$ and
$(f,f^{(k)})_0=\overline{\beta}_{k+1}$.
\item The polynomial $f^{(k)}$ is irreducible and its Puiseux
characteristic is
$\left(\frac{\beta_0}{l_k},\ldots,\frac{\beta_k}{l_k}\right)$ and
its semigroup is $\frac{\overline{\beta}_0}{l_k}\mathbf
N+\cdots+\frac{\overline{\beta}_k}{l_k}\mathbf N$.
\end{enumerate}
\end{Proposition}

\medskip

\noindent After \cite{Merle} and from the above Proposition, if
$$Q(f)=\langle q_1:m_1,\ldots,q_g:m_g \rangle$$
\noindent then
$$Q(f^{(k)})=\left \langle \frac{q_1}{l_k}:m_1,\ldots,
\frac{q_k}{l_k}:m_k \right \rangle.$$

\medskip

\noindent This inspires the next operation between the rational
Newton polygons.

\bigskip

\noindent {\bf The abrasion operation over rational Newton
polygons:}

\noindent Let ${\cal
N}=\sum_{i=1}^r\left\{\Teisssr{L_i}{M_i}{4}{2}\right \}$ be a
convenient rational Newton polygon with $r\geq 2$ and
$\frac{L_1}{M_1}<\cdots < \frac{L_r}{M_r}$. The {\em abrasion} of
${\cal N}$ is by definition
\begin{equation}
\label{abrasion polygon} {\cal A}({\cal
N}):=\sum_{i=1}^{r-1}\left\{\Teisssr{\tilde L_i}{ M_i}{4}{2}\right\}
\end{equation}
\noindent where

\medskip

\begin{equation}
\label{def3} \tilde
L_i:=\frac{1+M_1+\cdots+M_{r-1}}{1+M_1+\cdots+M_r}L_{i},
\end{equation}

\noindent for $1\leq i \leq r-1$.

\bigskip
\noindent Since $\frac{L_{i}}{M_{i}}<\frac{L_{i+1}}{M_{i+1}}$,
${\cal A}({\cal N})$  is again a convenient rational Newton polygon
written in its canonical form because $\frac{\tilde
L_i}{M_i}<\frac{\tilde L_{i+1}}{ M_{i+1}}$.

\medskip

\noindent The abrasion operation transforms a Newton polygon of
$r>1$ compact faces to a Newton polygon of $r-1$ compact faces.

\medskip

\noindent We denote by ${\cal A}^i({\cal N})$ the Newton polygon
obtained after applying $i$-times the abrasion operation to ${\cal
N}$. Observe that ${\cal A}^{i+1}({\cal N})={\cal A}^{i}({\cal
A}({\cal N}))={\cal A}({\cal A}^{i}({\cal N}))$. By convention we
put ${\cal A}^0({\cal N})={\cal N}$.

\medskip

\noindent More precisely if ${\cal
N}=\sum_{k=1}^r\left\{\Teisssr{L_k}{M_k}{4}{2}\right \}$ we denote
${\cal A}^i({\cal N})$ by $\sum_{k=1}^{r-i}\left\{\Teisssr{\tilde
L^{(i)}_k}{\tilde M^{(i)}_k}{6}{3}\right \}$.

\medskip

\begin{Properties}
Let $f$ be a branch with semigroup $\bigl\langle
\overline{\beta}_0,\overline{\beta}_1,\ldots
,\overline{\beta}_g\bigr\rangle$. Then
\begin{enumerate}
\item ${\cal A}({\cal N}_J(f))={\cal N}_J(f^{(g-1)})$,
\item ${\cal A}({\cal N}_J(f^{(i+1)}))={\cal N}_J(f^{(i)})$ for $1\leq i\leq g-1$
and
\item ${\cal A}^i({\cal N}_J(f))={\cal N}_J(f^{(g-i)})$ for $1\leq i\leq g-1$.
\end{enumerate}
\end{Properties}
\noindent \textbf{Proof.-} \noindent It is enough to prove that
${\cal A}({\cal N}_J(f))={\cal N}_J(f^{(g-1)})$. After \cite{Merle}
we have that ${\cal N}_J(f)=
\sum_{i=1}^g\left\{\Teisssr{L_i}{M_i}{4}{2}\right \}$, with
$L_i=(n_i-1)\overline{\beta}_i$ and $M_i=n_1\cdots n_{i-1}(n_i-1)$,
where $n_i=\frac{\hbox{\rm
g.c.d.}(\overline{\beta}_0,\ldots,\overline{\beta}_{i-1})}
{\hbox{\rm
g.c.d.}(\overline{\beta}_0,\ldots,\overline{\beta}_{i})}$. So ${\cal
A}({\cal N}_J(f))=\sum_{i=1}^{g-1}\left\{\Teisssr{\tilde
L_i}{M_i}{4}{2}\right \}$ where $\tilde L_i=\frac{1}{n_{g}}L_i$
since
$\frac{1+M_1+\cdots+M_{g-1}}{1+M_1+\cdots+M_g}=\frac{1}{n_{g}}$.

\medskip
\begin{Theorem}
\label{third main theorem} Let ${\cal
N}=\sum_{k=1}^r\left\{\Teisssr{L_k}{M_k}{4}{2}\right \}$ be a
special convenient integral Newton polygon satisfying

\begin{enumerate}
\item $1+\hbox{\rm ht}({\cal N})<\frac{L_1}{M_1}$,

\item ${\cal A}^i({\cal N})$ is a special convenient integral Newton polygon,
$\frac{\tilde L_1^{(i)}}{\tilde M_1^{(i)}} \in \mathbf N $ for all
$i\in \{0,\ldots,r-1\},$ and $(1+\tilde M_1^{(i)}+\tilde
M_2^{(i)}+\cdots+\tilde M_{r-i-1}^{(i)}) \frac{\tilde
L_{r-i}^{(i)}}{\tilde M_{r-i}^{(i)}}\in \mathbf N$ for all $i\in
\{0,\ldots,r-2\},$


\item $\hbox{\rm g.c.d.}\left(1+\hbox{\rm ht}({\cal A}^i({\cal N})), \tilde L_1^{(i)},
\ldots, \tilde L_{r-i}^{(i)})\right)=1$ for $0\leq i\leq r-1$,
\end{enumerate}

\noindent then there exists an irreducible $f\in\mathbf C\{x,y\}$
such that ${\cal N}={\cal N}_{J}(f)$, and its semigroup  is
$$(1+\hbox{\rm ht}({\cal N}))\mathbf N+\frac{L_1}{M_1}\mathbf N+
(1+M_1)\frac{L_2}{M_2}\mathbf N +\cdots +(1+M_1+M_2+\cdots+M_{r-1})
\frac{L_r}{M_r}\mathbf N.$$
\end{Theorem}

\noindent \textbf{Proof.-} \noindent We use induction on $r$.
According to Merle (\cite{Merle}) the Newton polygon ${\cal
N}=\left\{\Teisssr{L_1}{M_1}{4}{2}\right \}$ is the jacobian Newton
polygon of a branch $f(x,y)=0$ if and only if
$1+M_1<\frac{L_1}{M_1}$ and $\hbox{\rm
g.c.d.}(\frac{L_1}{M_1},1+M_1)=1$, which follows from the third
condition of Theorem~\ref{third main theorem} taking $r=1$ and
$i=0$.

\medskip

\noindent We suppose now that Theorem \ref{third main theorem} is
true for any special convenient integral Newton polygon with $r-1$
compact edges and let ${\cal
N}=\sum_{k=1}^r\left\{\Teisssr{L_k}{M_k}{4}{2}\right \}$ be a
special convenient integral Newton polygon satisfying the three
conditions of Theorem \ref{third main theorem}.  Then

\medskip
\noindent \noindent \textbf{Claim 1.} The abrasion operation ${\cal
A}({\cal N})$ of ${\cal N}$ also satisfies the hypothesis of 
Theorem~\ref{third main theorem}.

\medskip

\noindent The second and the third conditions are clearly satisfied.
Moreover by hypothe\-sis $1+\hbox{\rm ht}({\cal N})<\frac{L_1}{M_1}$
which is equivalent to $1+\hbox{\rm ht}({\cal A}({\cal
N}))<\frac{\tilde L_1}{M_1}$.

\bigskip
\noindent So there is an irreducible plane curve $\hat f$ such that
${\cal N}_J(\hat f)={\cal A}({\cal N})$.\\ Put $\gamma_0 \mathbf
N+\gamma_1 \mathbf N+\cdots+\gamma_{r-1}\mathbf N$ the semigroup of
$\hat f=0$, with $\gamma_i<\gamma_{i+1}$ for $0\leq i\leq r-2$.

\bigskip

\noindent We put $N=\frac{1+M_1+\cdots+M_{r}}{1+M_1+\cdots+M_{r-1}}$
and $\gamma_r=(1+M_1+M_2+\cdots+M_{r-1}) \frac{L_r}{M_r}$.

\medskip
\noindent \noindent \textbf{Claim 2.} The numbers $N$ and $\gamma_r$
are integers.

\medskip

\noindent The second condition of the Theorem gives, for $i=1$, that
${\cal A}^1({\cal N})$ is an integral Newton polygon, so in
particular $N$ divises $L_i$ for $i\in\{1,\ldots,r-1\}$ and the
third condition of the Theorem gives, for $i=1$,  $\hbox{\rm
g.c.d.}(1+M_1+\cdots+M_{r-1},\frac{L_1}{N},\ldots,\frac{L_{r-1}}{N})=1$
so $N=\hbox{\rm g.c.d.}(1+M_1+\cdots+M_{r},L_1,\ldots,L_{r-1})$.
Moreover $\gamma_r=(1+\tilde M_1^{(i)}+\tilde
M_2^{(i)}+\cdots+\tilde M_{r-i-1}^{(i)}) \frac{\tilde
L_{r-i}^{(i)}}{\tilde M_{r-i}^{(i)}}$ for $i=0$ so $\gamma_r \in
\mathbf N$ by hypothesis.

\medskip

\noindent \textbf{Claim 3.} The numbers
$N\gamma_0,N\gamma_1,\ldots,N\gamma_{r-1},\gamma_r$ form the minimal
set of generators of the semigroup of an irreducible plane curve.

\medskip

\noindent According to Bresinsky in \cite{Bresinsky} about the
characterization of the semigroups of plane branches (see the end of
Section \ref{section-Puiseux}), it is enough to prove the following
three conditions:

\begin{enumerate}
\item $\hbox{\rm g.c.d.}(N\gamma_0,N\gamma_1,\ldots,N\gamma_{r-1},\gamma_r)=1.$

\noindent We have $\hbox{\rm
g.c.d.}(N\gamma_0,N\gamma_1,\ldots,N\gamma_{r-1},\gamma_r)=\hbox{\rm
g.c.d.}(N,\gamma_r)$. Moreover $1+\hbox{\rm ht}({\cal N})$ is
divisible by $N$ and $N=\hbox{\rm
g.c.d.}(1+M_1+\cdots+M_{r},L_1,\ldots,L_{r-1})$ so $\hbox{\rm
g.c.d.}(N,L_r)=\hbox{\rm g.c.d.}(1+\hbox{\rm ht}({\cal
N}),L_1,\ldots,L_{r})=1$ and consequently $\hbox{\rm
g.c.d.}(N,\gamma_r)=1$ since $\gamma_r=\frac{L_r}{N-1}$.

\item For $1\leq i\leq r-1$, $\;\hbox{\rm g.c.d.}(N\gamma_0,N\gamma_1,
\ldots,N\gamma_{i-1}) > \hbox{\rm
g.c.d.}(N\gamma_0,N\gamma_1,\ldots,N\gamma_{i})$ by inductive
hypothesis and $\hbox{\rm g.c.d.}(N\gamma_0,N\gamma_1,
\ldots,N\gamma_{r-1})=N >1= \hbox{\rm
g.c.d.}(N\gamma_0,N\gamma_1,\ldots,N\gamma_{r-1},\gamma_r)$. This
proves the second condition of Bresinski.

\item For $1\leq i\leq r-2$, $\frac{\hbox{\rm
g.c.d.}(N\gamma_0,N\gamma_1,\ldots,N\gamma_{i-1})}{\hbox{\rm
g.c.d.}(N\gamma_0,N\gamma_1,\ldots,N\gamma_{i})}N\gamma_i<N\gamma_{i+1}$
by inductive hypothesis. To finish the proof of this third condition
observe that ${\cal N}_J(\hat f)={\cal A}({\cal N})$, so $\hbox{\rm
g.c.d.}(\gamma_0,\gamma_1,\ldots,\gamma_{r-2})=
\frac{1+M_1+\cdots+M_{r-1}}{1+M_1+\cdots+M_{r-2}}$ and

\begin{eqnarray*}
 \frac{\hbox{\rm
g.c.d.}(N\gamma_0,N\gamma_1,\ldots,N\gamma_{r-2})}{\hbox{\rm
g.c.d.}(N\gamma_0,N\gamma_1,\ldots,N\gamma_{r-1})}N\gamma_{r-1}&=&\hbox{\rm
g.c.d.}(\gamma_0,\gamma_1,\ldots,\gamma_{r-2})N\gamma_{r-1}\\
&=&(1+M_1+\cdots+M_{r-1}) \frac{L_{r-1}}{M_{r-1}}\\& < & \gamma_r.
\end{eqnarray*}

\noindent Finally note that the semigroup generated by
$N\gamma_0,\ldots,N\gamma_{r-1}, \gamma_r$ is exac\-tly

$$(1+\hbox{\rm
ht}({\cal N}))\mathbf N+\frac{L_1}{M_1}\mathbf N+
(1+M_1)\frac{L_2}{M_2}\mathbf N +\cdots +(1+M_1+M_2+\cdots+M_{r-1})
\frac{L_r}{M_r}\mathbf N.$$

\bigskip
\end{enumerate}

\begin{Corollary}\label{Cor3}
\label{characterization2} A special convenient integral Newton
polygon ${\cal N}=\sum_{k=1}^r\left\{\Teisssr{L_k}{M_k}{4}{2}\right
\}$ is the jacobian Newton polygon of a branch if and only if
verifies the next three conditions:

\begin{enumerate}
\item $1+\hbox{\rm ht}({\cal N})<\frac{L_1}{M_1}$,

\item ${\cal A}^i({\cal N})$ is a special convenient integral Newton polygon,
$\frac{\tilde L_1^{(i)}}{\tilde M_1^{(i)}} \in \mathbf N $ for all
$i\in \{0,\ldots,r-1\},$ and $(1+\tilde M_1^{(i)}+\tilde
M_2^{(i)}+\cdots+\tilde M_{r-i-1}^{(i)}) \frac{\tilde
L_{r-i}^{(i)}}{\tilde M_{r-i}^{(i)}}\in \mathbf N$ for all $i\in
\{0,\ldots,r-2\},$

\item $\hbox{\rm g.c.d.}\left(1+\hbox{\rm ht}({\cal A}^i({\cal N})), \tilde L_1^{(i)},
\ldots, \tilde L_{r-i}^{(i)})\right)=1$ for $0\leq i\leq r-1$.
\end{enumerate}

\end{Corollary}

\noindent {\bf Proof.}
It is enough to show that if ${\cal N}={\cal N}_J(f)$
for an irreducible distinguished Weierstrass polynomial $f(x,y)\in\Cn\{x\}[y]$
then ${\cal N}$ satisfies conditions 1--3.

\medskip

\noindent Let
$\bigl\langle\overline{\beta}_0,\overline{\beta}_1,\ldots,\overline{\beta}_r\bigr\rangle$
be the semigroup of $\{f=0\}$. Then by~(\ref{Me})
$${\cal N}=\sum_{k=1}^r\left\{\Teisssr{L_k}{M_k}{4}{2}\right\}=
\sum_{k=1}^r\left\{\Teisssr{(n_k-1)\overline{\beta}_k}
                                       {(n_k-1)n_1\dots n_{k-1}}{18}{9}
                                       \right\}\;.
$$

\noindent By above equality $L_1/M_1=\overline{\beta}_1$,
$(1+M_1+M_2+\cdots+M_{r-1})L_r/M_r=\overline{\beta}_r$  and
$\hbox{\rm g.c.d.}\left(1+\hbox{\rm ht}({\cal N}),
L_1,\ldots,L_{r}\right)= \hbox{\rm
g.c.d.}\left(\overline{\beta}_0,(n_1-1)\overline{\beta}_1,
\ldots,(n_r-1)\overline{\beta}_r\right)=1$ (see Property
\ref{arith_prop}). Hence conditions 2 and 3 are satisfied for
$i=0$.

\medskip

\noindent In order to check that conditions~2 and~3 are satisfied
for $i>1$ it is enough to observe that ${\cal A}({\cal N})={\cal
N}_J(f^{(r-1)})$ and apply an induction because an approximate root
$f^{(r-1)}$ is again an irreducible distinguished Weierstrass
polynomial.

\medskip
\begin{Corollary}
\label{Bre} Let ${\cal
N}=\sum_{k=1}^r\left\{\Teisssr{L_k}{M_k}{4}{2}\right \}$ be the
jacobian Newton polygon of a plane curve $f(x,y)=0$. Put
$\gamma_0:=1+\hbox{\rm ht}({\cal N})$,  and $\gamma_{i}:=(1+M_1+M_2+
\cdots+M_{i-1})\frac{L_i}{M_i}$ for $1\leq i \leq r$. Then $f$ is
irreducible if and only if the numbers
$\gamma_0,\gamma_1,\ldots,\gamma_r$ verify the three conditions of
Bresinsky. In such case $\gamma_0,\ldots,\gamma_r$ generate the
semigroup of $f$.
\end{Corollary}
\bigskip

\begin{Example}[Kuo's example]
Let $f(x,y)=(y^2-x^3)^2-x^7$ be Kuo's exam\-ple in \cite{Kuo}. Then
$\frac{\partial f}{\partial y}=4y(y^2-x^3)$, so $Q(f)=<6:1,7:2>$,
since $(f,y)_0=6$ and $(f,y^2-x^3)_0=14$. That is ${\cal
N}_J(f)=\{\Teis{6}{1}\}+\{\Teis{14}{2}\}$, so $1+\hbox{\rm ht}({\cal
N}_J(f))=4$, $\frac{L_1}{M_1}=6$ and $(1+M_1)\frac{L_2}{M_2}=14$.
Using Corollary \ref{Bre} we show that $f$ is not irreducible since
$\hbox{\rm g.c.d.}(4,6,14)\neq 1$. We can also show that $f$ is not
irreducible using the reduction operation since ${\cal R}({\cal
N}_J(f))=\{\Teis{8}{1}\}$ and the third condition, in Corollary
\ref{characterization}, is not true for $i=1$.
\end{Example}

\medskip

\begin{Remark}
In the case of a germ of plane irreducible curve (i.e. a branch),
Merle shows in \cite{Merle} that the datum of the jacobian Newton
polygon determines and is determined by the equisingularity class of
the curve. Note that if ${\cal N}$ is the jacobian Newton polygon of
a branch $f(x,y)=0$ of semigroup generated by
$\overline{\beta}_0,\ldots, \overline{\beta}_g$ then ${\cal N}$ has
$g$ compact faces $\Gamma_1,\ldots,\Gamma_g$ where the inclinations
of these faces form an increasing sequence. Note by ${P}_i$ the
segment which is the projection of $\Gamma_i$ on the vertical axis.
In particular the intersection point of $\Gamma_1$ and the vertical
axis is $(0,\overline{\beta}_0-1)$. If $l_i$ is the line containing
the point $(0,\overline{\beta}_0)$ and parallel to $\Gamma_i$,
${\cal B}_i=P_i\times \mathbf R$, and $r_i$(resp. $r_{i+1}$) is the
top line (resp. the bottom line) of ${\cal B}_i$ then the abscissa
of the intersection point of $l_i$ and $r_i$ is $\overline{\beta}_i$
and the abscissa of the intersection point of $l_i$ and $r_{i+1}$ is
$n_i\overline{\beta}_i$, where $n_i= \frac{\hbox{\rm
g.c.d.}(\overline{\beta}_0,\ldots,\overline{\beta}_{i-1})}
{\hbox{\rm
g.c.d.}(\overline{\beta}_0,\ldots,\overline{\beta}_{i})}$. Finally
note that the length of $P_i$ is exactly $n_1\ldots n_{i-1}(n_i-1)$.

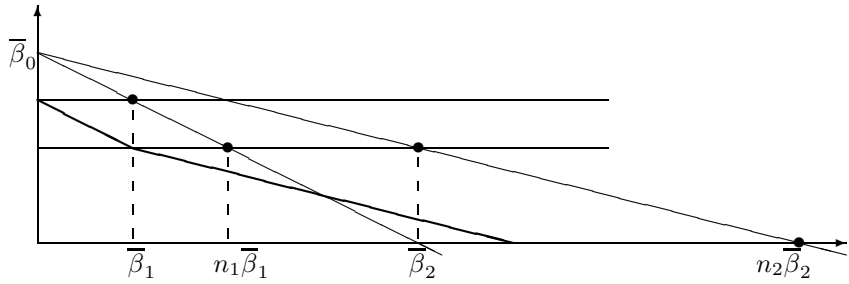
\setlength{\unitlength}{18pt}
\begin{picture}(20,7)(0,0)

\put(-0.6,3.8){$\overline\beta_0$}
\put(1.9,-0.6){$\overline\beta_1$}
\put(3.7,-0.6){$n_1\overline\beta_1$}
\put(7.8,-0.6){$\overline\beta_2$}
\put(15.1,-0.6){$n_2\overline\beta_2$}

\put(2,3){\makebox(0,0){$\bullet$}}
\put(4,2){\makebox(0,0){$\bullet$}}
\put(8,2){\makebox(0,0){$\bullet$}}
\put(16,0){\makebox(0,0){$\bullet$}}

\multiput(2,0.1)(0,0.5){6}{\line(0,1){0.2}}
\multiput(4,0.1)(0,0.5){4}{\line(0,1){0.2}}
\multiput(8,0.1)(0,0.5){4}{\line(0,1){0.2}}

\put(0,0){\vector(1,0){17}} \put(0,0){\vector(0,1){5}} \thicklines
\put(0,3){\line(2,-1){2}} \put(2,2){\line(4,-1){8}} \thinlines
\put(0,4){\line(2,-1){8.5}} \put(0,4){\line(4,-1){17}}
\put(0,3){\line(1,0){12}} \put(0,2){\line(1,0){12}}
\end{picture}

\end{Remark}

\bigskip

\noindent {\bf Conclusion.-} In this paper we characterize the
special convenient Newton polygons which are jacobian Newton
polygons of a branch. This allows us to give combinatorial criteria
of irreducibility of complex series in two variables. Their natures
are different from the classical criterion in \cite{Abhyankar}: we
draw the Newton polygon in the coordinates $(u,v)$ of the
discriminant $D(u,v)$ of the morphism defined by

\begin{equation}
\label{morphism} \begin{array}{rll} (x ,f)\colon (\mathbf
C^{2},0)&\longrightarrow& (\mathbf C^2,0)\\
(x,y)&\longrightarrow & (u,v):=(x,f(x,y)),
\end{array}
\end{equation}

\noindent and we check whether it verifies the geometrical
conditions of Corollary \ref{characterization} or Corollary
\ref{characterization2} or arithmetical conditions of Corollary
\ref{Bre}.

\medskip
\noindent By Corollary \ref{Cor1} we may assume that $f$ is a
Weierstrass polynomial. In this case it is not difficult to compute
an equation of the discriminant $D(u,v)$. We get this equation by
eliminating  $x$ and $y$ from
$$\left\{\begin{array}{r}
x=u\\
\\
f(x,y)=v\\
\\
\frac{\partial f(x,y)}{\partial y}=0
\end{array}\right.$$

\noindent and using the classical notion of resultant of two
polynomials in one variable we have:
\begin{eqnarray*}
D(u,v)&=& \hbox{\rm Result}_y\left(f(u,y)-v,\frac{\partial
f(u,y)}{\partial y}\right)\\&=& \hbox{\rm
Result}_y\left(f(u,y)-v,\frac{\partial (f(u,y)-v)}{\partial
y}\right)\\&=& \hbox{\rm Disc}_y(f(u,y)-v).
\end{eqnarray*}

\medskip

\noindent We think this approach provides effective methods of
checking the irreducibility of complex series in two variables.

\bigskip

\noindent

\noindent \textbf{Acknowledgments.} We would like to thank Andrzej
Lenarcik for his examples and comments and Arkadiusz P\l oski for
his comments and suggestions. Many thanks to Bernard Teissier, we
were also greatly helped in the development of Section
\ref{abrasion} by his questions and discussions.

\medskip
\noindent
{\small Evelia Rosa Garc\'{\i}a Barroso\\
Departamento de Matem\'atica Fundamental\\
Facultad de Matem\'aticas, Universidad de La Laguna\\
38271 La Laguna, Tenerife, Espa\~na\\
e-mail: ergarcia@ull.es}

\medskip

\noindent {\small Janusz Gwo\'zdziewicz\\
Department of Mathematics\\
Technical University \\
Al. 1000 L PP7\\
25-314 Kielce, Poland\\
e-mail: matjg@tu.kielce.pl}

\end{document}